\documentclass[12pt]{amsart}
\usepackage{amssymb,graphicx,color}
\newtheorem{thm}{Theorem}[section]
\newtheorem{lem}[thm]{Lemma}
\newtheorem{cor}[thm]{Corollary}
\theoremstyle{definition}
\newtheorem{dfn}{Definition}[section]
\newtheorem*{cnj}{Conjecture}
\theoremstyle{remark}
\newtheorem*{rem}{Remark}

\begin{document}
\title{Manifolds associated with $(\mathbb Z_2)^{n+1}$-colored regular graphs}
\author{Zhiqiang Bao}
\address{LMAM, School of Mathematical Sciences, Peking University, Beijing,
 100871,\newline P.R.China.}
\email{bzq@math.pku.edu.cn}
\author{Zhi L\"u}
\address{Institute of Mathematics, School of Mathematical Sciences, Fudan
 University,\newline Shanghai, 200433, P.R.China.}
\email{zlu@fudan.edu.cn}
\date{September 20, 2006}
\subjclass[2000]{Primary 57Q15, 57S17, 05C15, 05C25; Secondary 20F65, 55N91}
\keywords{Colored regular graph, skeletal expansion}
\thanks{The first author was supported by FANEDD Grant No. 200214 and NSFC
 Grant No. 10201002. The second author was supported by grants from NSFC (No.
 10371020 and No. 10671034) and JSPS (No. P02299).}
\begin{abstract}
In this article we describe a canonical way to expand a certain kind
of $(\mathbb Z_2)^{n+1}$-colored regular graphs into closed
$n$-manifolds by adding cells determined by the edge-colorings
inductively. We show that every closed combinatorial $n$-manifold
can be obtained in this way. When $n\leq 3$, we give simple
equivalent conditions for a colored graph to admit an expansion. In
addition, we show that if a $(\mathbb Z_2)^{n+1}$-colored regular
graph admits an $n$-skeletal expansion, then it is realizable as the
moment graph of an $(n+1)$-dimensional closed $(\mathbb
Z_2)^{n+1}$-manifold.
\end{abstract}
\maketitle

\section{Introduction}

In 1998, Goresky,  Kottwitz and MacPherson \cite{GKM} established
the GKM theory, which gives a direct link between equivariant
topology and combinatorics (or more precisely, a link between GKM
manifolds and GKM graphs). The most essential properties of a GKM
manifold $X$, such as the equivariant cohomology and Betti numbers
of $X$, can thus be explicitly expressed in terms of its
corresponding GKM graph. A lot of work on this subject has been
carried out since then (\cite{BGH}, \cite{GH}, \cite{GZ1},
\cite{GZ2}, \cite{GZ3}, \cite{GZ4}, \cite{MMP}, \cite{MP}). For
example, a series of papers by Guillemin and Zara showed that much
geometrical and topological information of a GKM manifold can be
read out from the corresponding GKM graph. In the GKM theory, the
action group is a torus. Alternately, when the action group is
chosen to be a mod 2-torus, there is an analogous theory, namely the
mod 2 GKM theory (see \cite{BGH} or \cite{L}), which has
successfully been applied to find the lower bound of the number of
fixed points and to study the equivariant cobordism classification
and the Smith problem in \cite{L}.

\vskip .2cm

The work of this article is mainly motivated by the GKM theory.
However, we shall carry out our work in the mod 2 category. Thus,
let us first recall some basic facts about the mod 2 GKM theory,
which inspired us to study colored regular graphs. Suppose that
$G=(\mathbb Z_2)^{n+1}$ (a mod 2-torus of rank $n+1$), which is also
an $(n+1)$-dimensional linear space over $\mathbb Z_2$. Let
$EG\rightarrow BG$ be the universal principal $G$-bundle, where
$BG=EG/G$ is the classifying space of $G$. It is well-known that
$BG=(\mathbb RP^\infty)^{n+1}$ and
$$H^*(BG;\mathbb Z_2)=\mathbb Z_2[t_0,t_1,\ldots,t_n]$$
is a $\mathbb Z_2$-polynomial ring over $t_0,t_1,\ldots,t_n$, where
the $t_i$'s are one-dimensional generators.

\vskip .2cm

Now suppose that $X$ is a closed manifold with an effective $G$-action such
that the fixed point set $X^G$ is finite, which implies $\dim X\geq n+1$ (see
\cite{AP}). Then the \emph{Borel construction} $EG\times_G X$ is defined as
the orbit space of the induced diagonal action $\Delta_G$ on $EG\times X$,
i.e.,
$$EG\times_G X=(EG\times X)/\Delta_G,$$
so that the equivariant cohomology of $X$ is
$$H^*_G(X;\mathbb Z_2)=H^*(EG\times_G X;\mathbb Z_2).$$
As shown in \cite{L}, to the $G$-manifold $X$, one can always assign
a regular graph $\Gamma$ with a label $\alpha(e)\in$
$\mathrm{Hom}(G,\mathbb Z_2)$ for each edge $e$ in $\Gamma$. The GKM
theory points out that when the $G$-manifold $X$ satisfies certain
conditions, its equivariant cohomology is a free module over
$H^*(BG;\mathbb Z_2)$ and can be explicitly read out from $(\Gamma,
\alpha)$:
$$H^*_G(X;\mathbb Z_2)=\{f:V_\Gamma\rightarrow\mathbb
 Z_2[t_0,...,t_n]~|~f(v)\equiv f(w)\mod\alpha(e)~\text{ for }e\in E_v\cap
 E_w\}$$
where $V_\Gamma$ denotes the vertex set of $\Gamma$ and $E_v$
denotes the set of all edges in $\Gamma$ that are adjacent to $v$.

\vskip .2cm

In the extreme case where $\dim X=n+1$, the associated graph
$\Gamma$ with $X^G$ as its vertex set can be derived as follows.

\vskip .2cm

Every real irreducible representation of $G$ is one-dimensional, so it must be
a homomorphism $\rho:G\rightarrow GL(1,\mathbb R)$. Moreover, for each $g\in
G$, $\rho(g)=(\pm 1)$ since $g$ is an involution. Therefore, by identifying
the multiplicative group $\{(1),(-1)\}$ with $\mathbb Z_2$, one can identify
the set of all real irreducible representations of $G$ with
$\mathrm{Hom}(G,\mathbb Z_2)$. Notice that $\mathrm{Hom}(G,\mathbb Z_2)\cong
G$ is also an $(n+1)$-dimensional linear space over $\mathbb Z_2$.

\vskip .2cm

Now choose a point $v\in X^G$, and choose a linear identification
$\phi_v: T_v X\rightarrow\mathbb R^{n+1}$.  Then the induced tangent
representation at $v$ can be expressed as an $(n+1)$-dimensional
real representation $\rho:G\rightarrow GL(n+1,\mathbb R)$, which can
be further split into $n+1$ real irreducible representations
$$\rho_0,\rho_1,...,\rho_n\in \mathrm{Hom}(G,\mathbb Z_2).$$
Clearly the set $\{\rho_0,\rho_1,\ldots,\rho_n\}$ is independent of
the choice of $\phi_v$. Since $G$ acts effectively on $X$,
$\rho_0,\rho_1,\ldots,\rho_n$ are linearly independent in
$\mathrm{Hom}(G,\mathbb Z_2)$, so they form a basis of
$\mathrm{Hom}(G,\mathbb Z_2)$.

\vskip .2cm

For each $\rho_i$, the subgroup $\mathfrak g_i=\ker\rho_i$ is of
codimension one in $G$, and the $1$-dimensional component of the
fixed set of $\mathfrak g_i$ acting on $X$ is a circle $S_i$ (c.f.
\cite{AP} or \cite{CF}). Since $G/\mathfrak g_i\cong\mathbb Z_2$,
the $G$-action on $X$ induces an involution on $S_i$, so $S_i$
contains only two fixed points $v$ and $w$ in $X^G$. Let us label
$S_i$ by $\rho_i$ and choose an orientation on $S_i$.  Then one gets
two \emph{oriented} edges with endpoints $v, w$ and with the same
``label'' $\rho_i$. Since the labeled oriented edges always appear
in pairs with reversing orientations, one can merge each pair into
an \emph{unoriented} edge, so that each $S_i$ determines one edge
with two $G$-fixed points as its endpoints (i.e., a $1$-valent
graph). Then, our desired graph $\Gamma$ is the union of all these
$1$-valent graphs such that the vertex set of $\Gamma$ consists of
all fixed points of $X^G$. $\Gamma$ is called the \emph{moment
graph} and the ``labelling map''
$\alpha:E_\Gamma\rightarrow\mathrm{Hom}(G,\mathbb Z_2)$ is called
the \emph{axial function}, where $E_\Gamma$ denotes the set of all
unoriented edges.

\vskip .2cm

It is obvious that in this simple case, the moment graph $\Gamma$ is
a regular $(n+1)$-valent graph, and the axial function
$\alpha:E_\Gamma\rightarrow \mathrm{Hom}(G,\mathbb Z_2)$ satisfies
the following two conditions:

\begin{enumerate}
\item{for each vertex $v$ in $\Gamma$, if $e_0,\ldots,e_n$ are all edges
 adjacent to $v$, then $\alpha(e_0),\ldots,\alpha(e_n)$ linearly span
 $\mathrm{Hom}(G,\mathbb Z_2)$, so they are all nontrivial and distinct};
\item{for each edge $e$ in $\Gamma$, let $\mathfrak g_e=\ker\alpha(e)$ and
 $v,w$ denote the two endpoints of $e$. If $\{e_0,\ldots,e_n\}$ and
 $\{e'_0,\ldots,e'_n\}$ are the sets of edges adjacent to $v$ and $w$ respectively,
 then
 $$\{\alpha(e_0)|_{\mathfrak g_e},\ldots,\alpha(e_n)|_{\mathfrak
  g_e}\}=\{\alpha(e'_0)|_{\mathfrak g_e},\ldots,\alpha(e'_n)|_{\mathfrak
  g_e}\}$$
 (because they form irreducible decompositions of two equivalent
 tangent representations of $\mathfrak g_e$ at $v$ and $w$ respectively). It is
 not difficult to verify that this is equivalent to
 \begin{eqnarray*}
  &&\{\mathrm{Span}(\alpha(e_0),\alpha(e)),\ldots,
   \mathrm{Span}(\alpha(e_n),\alpha(e))\}\\
  &=&\{\mathrm{Span}(\alpha(e'_0),\alpha(e)),\ldots,
   \mathrm{Span}(\alpha(e'_n),\alpha(e))\}.
 \end{eqnarray*}}
\end{enumerate}

\vskip .2cm

Conversely, given a finite connected regular graph $\Gamma$ together
with some axial function $\alpha$ satisfying the previous two
conditions, a natural question is: \emph{what information on the
topology of manifolds can we extract from the pair
$(\Gamma,\alpha)$}? This seems to be an open-ended question.
Recently, some interesting and beautiful work on this subject has
been done by some mathematicians. For example, in \cite{GZ2}
Guillemin and Zara obtained some purely combinatorial analogues of
the main results in the GKM theory. In addition, they also obtained
a realization theorem for abstract GKM-graphs: for an abstract GKM
graph $(\Gamma,\alpha)$, there exists a complex manifold $X$ and a
GKM action of the torus $T$ on $X$ such that $(\Gamma,\alpha)$ is
its GKM graph. Note that the constructed complex manifold $X$ is not
compact, not equivariantly formal, and admits no canonical
compactification.

\vskip .2cm

This article will deal with the problem from a different viewpoint,
in the mod $2$ category. We shall introduce a ``skeletal expansion
technique'' for any given abstract $G$-colored finite connected
regular $(n+1)$-valent graphs $(\Gamma, \alpha)$, where $\alpha$
satisfies the above two conditions. Specifically, we shall give a
canonical way to inductively attach cells on $\Gamma$ to form a cell
complex $K$ such that the $1$-skeleton of $K$ is just $\Gamma$. We
will call this technique the \emph{skeletal expansion}, and we shall
use it to carry out our work as follows.

\vskip .2cm

First, we determine under what conditions this procedure of
cell-gluing can be performed to the end until one obtains a closed
manifold (see Theorem~\ref{th:manifold}). We further show that every
closed combinatorial $n$-manifold can be obtained in this way (see
Theorem~\ref{th:colorable}).

\vskip .2cm

Next, we consider the realization problem: can we always construct a
$G$-action on an $(n+1)$-dimensional closed manifold such that its
moment graph is exactly the given $G$-colored regular $(n+1)$-valent
graph $(\Gamma, \alpha)$? To do this, our strategy is first to try
to construct the orbit space of an action from $(\Gamma, \alpha)$.
Then by the work of Davis and Januszkiewicz \cite{DJ} on the
reconstruction of small covers (which is stated briefly below), we
shall use the obtained orbit space to construct the desired
$G$-action. We shall show that under certain conditions, any
abstract $G$-colored finite connected regular $(n+1)$-valent graph
$(\Gamma, \alpha)$ is realizable as the moment graph of an
$(n+1)$-dimensional closed $G$-manifold (see
Theorem~\ref{realization}).

\vskip .2cm

An \emph{$(n+1)$-dimensional convex polytope} $P\subset\mathbb
R^{n+1}$ is just an $(n+1)$-dimensional compact manifold which is
the intersection of a finite set of half-spaces in $\mathbb
R^{n+1}$. It naturally becomes a cell decomposition of $B^{n+1}$ and
has a minimum set of defining half-spaces, each of which corresponds
to an $n$-dimensional face of $P$. Denote the set of all these
codimension one faces by $\mathcal F(P)$. If every vertex is
surrounded by exactly $n+1$ faces in $\mathcal F(P)$, then $P$ is
called a \emph{simple} convex polytope. Clearly, the $1$-skeleton of
each $(n+1)$-dimensional simple convex polytope is a regular
$(n+1)$-valent graph $\Gamma(P)$.

\vskip .2cm

If $P$ is an $(n+1)$-dimensional simple convex polytope, a
\emph{characteristic function} is a map $\lambda:\mathcal F(P)\to
G=(\mathbb Z_2)^{n+1}$ such that the $n+1$ faces in $\mathcal F(P)$
adjacent to each vertex are sent to $n+1$ linearly independent
vectors in $G$. Each such $\lambda:\mathcal F(P)\to G$ is dual to an
axial function
$\alpha_\lambda:E_{\Gamma(P)}\to\mathrm{Hom}(G,\mathbb Z_2)$, such
that for each edge $e\in E_{\Gamma(P)}$ and those faces
$F\in\mathcal F(P)$ containing $e$,
$\alpha_\lambda(e)(\lambda(F))=0$. It is easy to see that both
$\lambda$ and $\alpha_\lambda$ are determined by each other.

\vskip .2cm

An $(n+1)$-dimensional closed manifold $X$ is called a \emph{small
cover} over $P$ if it admits an effective and locally standard
$G$-action such that the orbit space $X/G=P$. In \cite{DJ}, Davis
and Januszkiewicz observed that for each simple convex polytope $P$
and each characteristic function $\lambda:\mathcal F(P)\to G$, there
is a canonical way to construct a small cover $X(\lambda)$ over $P$,
and every small cover over $P$ can be reconstructed in this way.
Moreover, $\Gamma(P)$ is exactly the moment graph for $X(\lambda)$
and $\alpha_\lambda$ is the corresponding axial function.

\vskip .2cm

In the case of small covers, reconstructing the orbit spaces and
$G$-manifolds from moment graphs is simple. In fact, the orbit space
$X/G$ of a small cover $X$ is bounded by $S^n$ with a very nice
regular cell decomposition corresponding to $\partial P$, and this
cell complex is topologically dual to a triangulation of $S^n$. The
skeletal expansion of the moment graph $(\Gamma,\alpha)$ of $X$ will
exactly reproduce the pre-defined cell decomposition on $\partial
P$, so that this can lead us to recover the orbit space $P$ and its
characteristic function $\lambda$, and thus $X$ can be
reconstructed. More generally, if $X$ is not a small cover, the
problem of reconstructing the $G$-manifold $X$ and its orbit space
compatible with the moment graph becomes much more complicated. We
remark that if the orbit space $X/G$ is a compact $(n+1)$-manifold
with boundary such that the pre-image of each component of the
boundary of $X/G$ contains a fixed point, then the skeletal
expansion of the moment graph $(\Gamma,\alpha)$ will produce a cell
decomposition of the boundary of $X/G$.

\vskip .2cm

\textbf{Acknowledgements.} The authors would like to thank M.
Masuda, F. Luo and X. R. Zhang for very helpful discussions leading
to the results in dimensions two and three. The authors also would
like to thank the referee, who gave many valuable suggestions and
comments.

\section{Notation and main results}

\begin{figure}[ht]
 \includegraphics{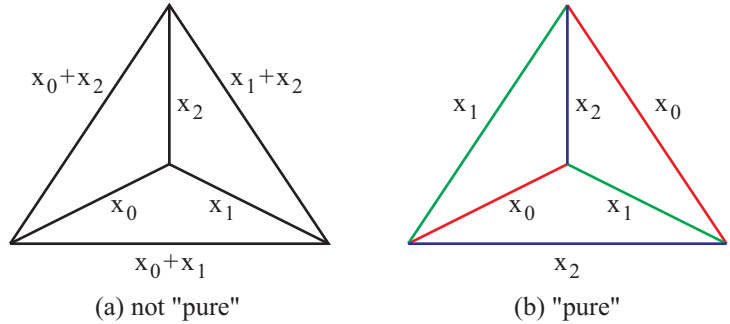}
 \caption{Examples of $G$-colored graphs}\label{fg:coloring}
\end{figure}

In this section, we will first define some standard notation and
state our main results. Suppose $G=(\mathbb Z_2)^{n+1}$ and $\Gamma$
is a finite connected regular $(n+1)$-valent graph. A
\emph{$G$-coloring} is a map
$$\alpha:E_\Gamma\to\mathrm{Hom}(G,\mathbb Z_2),$$
such that at each vertex, the colors of the $n+1$ edges adjacent to
it are linearly independent vectors (here \emph{the color of an edge
$e$} is simply the vector $\alpha(e)$). Notice that for some
technical reason, the edge colors are in $\mathrm{Hom}(G,\mathbb
Z_2)\cong G$, but \emph{not} in $G$. If in addition the total image
of $\alpha$ contains only $n+1$ vectors $x_0,\ldots,x_n$, then we
call $\alpha$ a \emph{``pure'' $G$-coloring}, which is just the
\emph{ordinary edge coloring} (see Fig.~\ref{fg:coloring}).

For each subgraph $\Delta\subseteq\Gamma$, let $\mathrm{Span}~\alpha(\Delta)$
denote the linear space spanned by all colors of edges in $\Delta$, and call
it \emph{the color of $\Delta$}. Let
$$\dim_\alpha\Delta=\dim\mathrm{Span}~\alpha(\Delta).$$

If $K$ is a triangulation of a manifold $M$ and the graph $\Gamma$
is its $1$-skeleton, then each $k$-dimensional simplex $\Delta\in K$
can be characterized by the $1$-skeleton of $\partial\Delta$, which
is a subgraph of $\Gamma$. Motivated by this idea, we introduce the
following concept of $k$-nest for colored graph.

\begin{dfn}
We say that a proper subgraph $\Delta\subsetneqq\Gamma$ is a
\emph{$k$-dimensional colored nest} (or just a \emph{$k$-nest})
associated with $(\Gamma,\alpha)$ if
\begin{enumerate}
 \item{$\Delta$ is connected;}
 \item{$\dim_\alpha\Delta=k$;}
 \item{for every connected subgraph $\Delta'\supsetneqq\Delta$,
 $\dim_\alpha\Delta'>k$.}
\end{enumerate}
If $\Delta'\supseteq\Delta$, and they are both colored nests, then
we say that $\Delta$ is a \emph{face} of $\Delta'$. Let
$K(\Gamma,\alpha)$ denote the set of all colored nests of
$(\Gamma,\alpha)$. We call
$$K_k(\Gamma,\alpha)=\{\Delta\in K(\Gamma,\alpha)~|\dim_\alpha\Delta\leq k\}$$
the \emph{$k$-dimensional colored nest-skeleton} (or just the
\emph{$k$-skeleton}) of $K(\Gamma,\alpha)$.
\end{dfn}

\begin{rem}
Suppose $\mathrm{Hom}(G,\mathbb Z_2)$ is spanned by $n+1$ vectors
$x_0,\ldots,x_n$. Let $\mathcal S(G,\mathbb Z_2)$ denote the $\mathbb
Z_2$-exterior algebra generated by these vectors, which consists of all
square-free (i.e., we require $x_0^2=0,\ldots,x_n^2=0$ in $\mathcal
S(G,\mathbb Z_2)$) $\mathbb Z_2$-polynomials of $x_0,\ldots,x_n$, and let
$\mathcal S^k(G,\mathbb Z_2)$ denote its linear subspace spanned by all degree
$k$ monomials, namely
$$x_0\cdots x_{k-1},x_0\cdots x_{k-2}x_k,\ldots,x_{n-k+1}\cdots x_n.$$
Then $\mathcal S^1(G,\mathbb Z_2)=\mathrm{Hom}(G,\mathbb Z_2)$, so
the color of each edge can also be considered as an element in
$\mathcal S(G,\mathbb Z_2)$. Furthermore, for any $m$-dimensional
linear subspace $L\subseteq\mathcal S^1(G,\mathbb Z_2)$, choose a
basis $y_1,\ldots,y_m$ for $L$, one can easily verify that
$\eta_L=y_1\cdots y_m\in\mathcal S^m(G,\mathbb Z_2)$ is independent
of the choice of $y_1,\ldots,y_m$. Therefore, for each connected
subgraph $\Delta\subseteq\Gamma$, $L=\mathrm{Span}~\alpha(\Delta)$
determines a unique homogeneous $\mathbb Z_2$-polynomial
$\eta_L\in\mathcal S^{\dim_\alpha\Delta}(G,\mathbb Z_2)$. In
particular, one can extend the coloring $\alpha$ to a ``coloring''
$$\tilde\alpha:K(\Gamma,\alpha)\to\mathcal S(G,\mathbb Z_2).$$
\end{rem}

Clearly a $0$-nest $v$ is just a vertex of $\Gamma$. A $1$-nest $e$ is an edge
in $E_\Gamma$ together with its two ends. A $2$-nest $\gamma$ is a maximal arc
or maximal circle such that $\dim_\alpha\gamma=2$.

\begin{dfn}
If for each edge $e_1$  with endpoints $v,w$ and each edge $e_0$
adjacent to $v$, there exists one and only one edge $e_2$ adjacent
to $w$ such that
$$\mathrm{Span}(\alpha(e_0),\alpha(e_1))
 =\mathrm{Span}(\alpha(e_1),\alpha(e_2)),$$
then we call $\alpha$ a \emph{``good'' coloring}. Clearly ``pure''
colorings are always ``good'' colorings.
\end{dfn}

\begin{rem}
Notice that the $\alpha(e_i)$'s can be considered as elements in
$\mathcal S^1(G,\mathbb Z_2)$. Assume that the colors of edges
adjacent to each vertex are linearly independent, then the above
equation is equivalent to
$\alpha(e_0)\alpha(e_1)=\alpha(e_1)\alpha(e_2)$, and is also
equivalent to
$$\prod_{e\in E_v\backslash\{e_1\}}\alpha(e)\equiv
 \prod_{e\in E_w\backslash\{e_1\}}\alpha(e)\mod\alpha(e_1),$$
where $E_p$ denotes the set of all edges in $\Gamma$ adjacent to a
vertex $p$.  This also implies that there must be a unique bijection
$\theta_{e_1}: E_v\longrightarrow E_w$ such that for any $e\in
E_v\backslash\{e_1\}$,
$$\alpha(e)\equiv\alpha(\theta_{e_1}(e))\mod \alpha(e_1).$$
The collection $\theta=\{\theta_e\big\vert e\in E_\Gamma\}$ is
called a \emph{connection} of $(\Gamma, \alpha)$. Notice also that
if $X$ is an $(n+1)$-dimensional closed manifold with an effective
and locally standard $G$-action such that the fixed point set $X^G$
is finite, then the orbit space $X/G$ is a nice manifold with
corners (see, \cite{D} and \cite{LM}) and each axial function on its
moment graph is in fact a ``good'' coloring.
\end{rem}

\begin{lem}\label{lm:regular}
Suppose $n\geq 2$ and $(\Gamma,\alpha)$ is a finite connected
$G$-colored regular graph. Then $\alpha$ is ``good'' if and only if
each $k$-nest of $K(\Gamma,\alpha)$ is a connected regular
$k$-valent subgraph.
\end{lem}

\begin{proof}
Suppose $\alpha$ is ``good''. Fix a $k$-nest $\Delta$. For any
vertex $p$, let $E(p)$ denote the set of all edges \emph{inside}
$\Delta$ adjacent to $p$, and let $L_p$ denote the linear space
spanned by $\{\alpha(e)~|~e\in E(p)\}$. Now let us show that for any
two vertices $u$, $v$ connected by an edge $e_1\in\Delta$,
$L_u=L_v$.

If this is not true, suppose without loss of generality that there
is an edge $e_0\in E(u)$ such that $\alpha(e_0)\not\in L_v$. Then
the connection of $\alpha$ points out an edge $e_2\in E_\Gamma$
adjacent to $v$ such that
$$\mathrm{Span}(\alpha(e_0),\alpha(e_1))=\mathrm{Span}(\alpha(e_1),\alpha(e_2)).$$
Since
$\alpha(e_2)\in\mathrm{Span}(\alpha(e_0),\alpha(e_1))\subseteq\mathrm{Span}~\alpha(\Delta)$
and $\Delta$ is maximal, $e_2$ must be an edge inside $\Delta$, so
$e_2\in E(v)$ and $\alpha(e_2)\in L_v$, which implies
$$\alpha(e_0)\in\mathrm{Span}(\alpha(e_1),\alpha(e_2))\subseteq L_v,$$
a contradiction.

Since $\Delta$ is connected, this implies that $L_p$ is the same for
all vertices $p$ in $\Delta$. Therefore for any vertex $p$ in
$\Delta$, $L_p=\mathrm{Span}~\alpha(\Delta)$. However, this is a
 $k$-dimensional linear space, thus the
valence of $p$ in $\Delta$ must be always $k$. Notice that for
$0<l\leq n+1$, the colors on $l$ edges emerging from the same vertex
are always linearly independent so they  span a $l$-dimensional
linear space.

On the other hand, suppose every $2$-nest is a regular $2$-valent graph. For
each edge $e_1$ with endpoints $v,w$ and each edge $e_0$ adjacent to $v$, let
$\tilde\gamma$ be the graph formed by edges with colors inside
$\mathrm{Span}(\alpha(e_0),\alpha(e_1))$ together with their vertices, and let
$\gamma$ be its connected component containing $e_0$ and $e_1$. Then $\gamma$
is a $2$-nest, so the valence of $\gamma$ at $w$ equals two. This implies that
there is exactly one edge $e_2\neq e_1$ adjacent to $w$,  such that
$\alpha(e_2)\in\mathrm{Span}(\alpha(e_0),\alpha(e_1))$. Hence $\alpha$ is
``good''.
\end{proof}

\begin{rem}
The proof of Lemma~\ref{lm:regular} actually shows that if $\alpha$
is ``good'', then for any linear subspace $L$, the edges in
$\alpha^{-1}(L)=\{e\in E_\Gamma|\alpha(e)\in L\}$ together with
their vertices also forms a disjoint union of some regular graphs.
\end{rem}

Recall that a \emph{regular cell decomposition} of a topological space $M$ is
a cell decomposition $K$ such that
\begin{enumerate}
 \item{for each cell $\Delta\in K$, let $B$ be the unit \emph{open} ball with
  dimension equal to $\dim\Delta$, and let $S$ be its boundary sphere, then
  the characteristic map $\phi:B\to\Delta\subseteq M$ extends to a
  homeomorphism $\phi:(\overline B,S)\to(\overline\Delta,\partial\Delta)$
  where $\partial\Delta=\overline\Delta\setminus\Delta\subseteq M$;}
 \item{for each cell $\Delta\in K$, its boundary sphere $\partial\Delta$ is
  the union of some cells in $K$.}
\end{enumerate}

\begin{dfn}
Suppose $k\leq n$ and $M$ is a topological space which has a
$k$-dimensional finite regular cell decomposition $K$. If there is a
one-to-one correspondence $\kappa:K\to K_k(\Gamma,\alpha)$
preserving dimensions and face relations, then $(M,K)$ is called a
\emph{$k$-skeletal expansion of $(\Gamma,\alpha)$}, and
$K_k(\Gamma,\alpha)$ is called a \emph{$G$-colored cell
decomposition of $M$}. If $k=n$ and $M$ is a closed manifold, then
$M$ is called a \emph{$(\mathbb Z_2)^{n+1}$-colorable manifold}.
\end{dfn}

\begin{rem}
If $\Gamma$ is the $1$-skeleton of an $(n+1)$-dimensional simple
convex polytope $P$, and $\alpha:E_\Gamma\to\mathrm{Hom}(G,\mathbb
Z_2)$ is dual to a characteristic function of $P$ (c.f. section 1),
then $\partial P$ with the induced cell decomposition is an
$n$-skeletal expansion of $(\Gamma,\alpha)$.
\end{rem}

Obviously, $1$-skeletal expansions always exist. One can simply
choose the vertices of $\Gamma$ as the $0$-cells and choose the
edges in $E_\Gamma$ as the $1$-cells.

\begin{lem}\label{lm:skeletal}
Suppose $n\geq 2$ and the coloring $\alpha$ is ``good''. Then
$(\Gamma,\alpha)$ always has a $2$-skeletal expansion.
\end{lem}

\begin{proof}
By Lemma \ref{lm:regular}, each $2$-nest is a connected regular
$2$-valent graph. Therefore all 2-nests must be embedded circles.
Furthermore, one can glue discs to each of these circles in the
graph to get a $2$-skeletal expansion.
\end{proof}

\begin{rem}
If $n\geq 2$ and $\alpha:E_\Gamma\to\mathrm{Hom}(G,\mathbb Z_2)$ is
``pure'', then $(\Gamma,\alpha)$ always has a $2$-skeletal
expansion. Notice that a $2$-nest in this case must be a simple
closed curve in $\Gamma$ consisting of an even number of edges with
alternating edge-coloring. We call them \emph{bi-colored circles}.
\end{rem}

 Higher dimensional skeletal expansions can be
built up inductively using the following lemma, which can be easily
proved similarly.

\begin{lem}\label{lm:induction}
Suppose that $\alpha$ is ``good'', $(\Gamma,\alpha)$ has a
$k$-skeletal expansion $(M,K)$, and $\kappa:K\to K_k(\Gamma,\alpha)$
is a one-to-one correspondence. For each $k+1$ edges
$e_0,\ldots,e_k\in E_\Gamma$ sharing one common vertex $v$, let
$\Delta$ be the unique $(k+1)$-nest containing $e_0,\ldots,e_k$, and
let
$$F^k(e_0,\ldots,e_k)=\bigcup_{\kappa(\sigma)\subseteq\Delta}\sigma.$$
If each such $F^k(e_0,\ldots,e_k)$ is a $k$-sphere, then one can
glue a $(k+1)$-cell to each $k$-sphere and obtain a $(k+1)$-skeletal
expansion of $(\Gamma,\alpha)$.\hfill\qed
\end{lem}

In general, $M$ is only a topological space but not a manifold, so it would be
quite natural to ask the following three questions:
\begin{enumerate}
 \item[(Q1)]{When will a finite connected regular $(n+1)$-valent graph
  $\Gamma$ with a ``good'' $(\mathbb Z_2)^{n+1}$-coloring admit an
  $n$-skeletal expansion?}
 \item[(Q2)]{If $(\Gamma,\alpha)$ admits an $n$-skeletal expansion $(M,K)$,
  when will $M$ be a closed manifold?}
 \item[(Q3)]{What kind of $n$-dimensional manifold can be $(\mathbb
  Z_2)^{n+1}$-colorable?}
\end{enumerate}

In the next two sections, we will give a complete answer to these
questions in dimensions $2$ and $3$. For dimensions greater than
$3$, question (Q1) is still open, due to the difficulties in
recognizing regular cell decompositions of $S^3$. However, questions
(Q2) and (Q3) are solved completely by the following two theorems.

\begin{thm}\label{th:manifold}
Suppose $G=(\mathbb Z_2)^{n+1}$, $(\Gamma,\alpha)$ is a $G$-colored
finite connected regular $(n+1)$-valent graph in which $\alpha$ is a
``good'' coloring, and $(M,K)$ is an $n$-skeletal expansion of
$(\Gamma,\alpha)$. Then $M$ is an $n$-dimensional closed manifold.
\end{thm}

An $n$-dimensional \emph{closed combinatorial manifold} is a connected
topological space with a finite $n$-dimensional regular cell decomposition,
such that the link complex of each $k$-cell is an embedded $(n-k-1)$-sphere.
Notice that every closed differentiable manifold is a closed combinatorial
manifold (see e.g. \cite{W}), and every closed combinatorial manifold is a
closed (topological) manifold.

\begin{thm}\label{th:colorable}
Every $n$-dimensional closed combinatorial manifold is $(\mathbb
Z_2)^{n+1}$-colorable. Moreover, one can choose ``pure'' colorings
for doing this.
\end{thm}

\begin{rem}
It is not difficult to see that this is also a necessary condition for being
colorable, namely an $n$-dimensional closed manifold is $(\mathbb
Z_2)^{n+1}$-colorable \emph{if and only if} it is a closed combinatorial
manifold.
\end{rem}

Although a manifold constructed as above does not naturally carry
any group actions, in some sense, it is associated with the boundary
of the orbit space of a manifold with group actions. Thus we can
pose the following realization problem:

\begin{enumerate}
 \item[(Q4)]{ Let $(\Gamma,\alpha)$ be a $(\mathbb Z_2)^{n+1}$-colored finite
 connected regular $(n+1)$-valent graph such that  $\alpha$ is a
 ``good'' coloring. Under what conditions will $(\Gamma,\alpha)$ be realizable
 as the moment graph of a $(\mathbb Z_2)^{n+1}$-action on an $(n+1)$-dimensional
 closed manifold?}
\end{enumerate}

We know by Theorem~\ref{th:manifold} that if $(\Gamma,\alpha)$
admits an $n$-skeletal expansion $(M,K)$, then $M$ is an
$n$-dimensional closed manifold. By a simple trick, we can modify
$M$ into the boundary of another manifold $N$ and turn $N$ into the
orbit space of a $G$-manifold. This leads us to the following
result.

\begin{thm}\label{realization}
Suppose $(\Gamma,\alpha)$ is a $(\mathbb Z_2)^{n+1}$-colored finite
connected regular $(n+1)$-valent graph such that  $\alpha$ is a
``good'' coloring.  If $(\Gamma,\alpha)$ admits  an $n$-skeletal
expansion  $(M,K)$, then it is realizable as the moment graph of
some $(\mathbb Z_2)^{n+1}$-action on an $(n+1)$-dimensional closed
manifold.
\end{thm}

Although we have not achieved a complete answer for the realization
problem, it is very tempting to make the following conjecture:

\begin{cnj}
For each $(\mathbb Z_2)^{n+1}$-colored finite connected regular
$(n+1)$-valent graph $(\Gamma,\alpha)$ with $\alpha$ being a
``good'' coloring, $(\Gamma,\alpha)$ can be realized as the moment
graph of some $(\mathbb Z_2)^{n+1}$-action on an $(n+1)$-dimensional
closed manifold.
\end{cnj}

\section{$(\mathbb Z_2)^3$-colorable closed surfaces}

Let us first examine the case $n=2$, which is surprisingly simple.

\begin{thm}\label{th:surface}
Suppose $G=(\mathbb Z_2)^3$, $\Gamma$ is a finite connected regular $3$-valent
graph, and $\alpha:E_\Gamma\to\mathrm{Hom}(G,\mathbb Z_2)$ is a ``good''
coloring. Then $(\Gamma,\alpha)$ has a $2$-skeletal expansion $(M,K)$ in which
$M$ is a closed surface.
\end{thm}

\begin{proof}
By Lemma \ref{lm:skeletal}, $(\Gamma,\alpha)$ must have a
$2$-skeletal expansion $(M,K)$. Moreover, since each edge is used
exactly twice by $2$-cells, this gives a pairing of arcs on the
boundaries of those $2$-cells. Therefore, gluing these cells
together along their boundaries will give a closed surface.
\end{proof}

As an example of application, the two graphs in
Fig.~\ref{fg:coloring} both satisfy the above conditions. It is not
difficult to see that the surface corresponding to
Fig.~\ref{fg:coloring} (a) is a 2-sphere, while the surface for
Fig.~\ref{fg:coloring} (b) is a real projective plane $P^2$.

\begin{thm}
Every closed surface $M$ is $(\mathbb Z_2)^3$-colorable. Moreover, one can
choose planar graphs and ``pure'' colorings for doing this.
\end{thm}

\begin{figure}[ht]
 \includegraphics{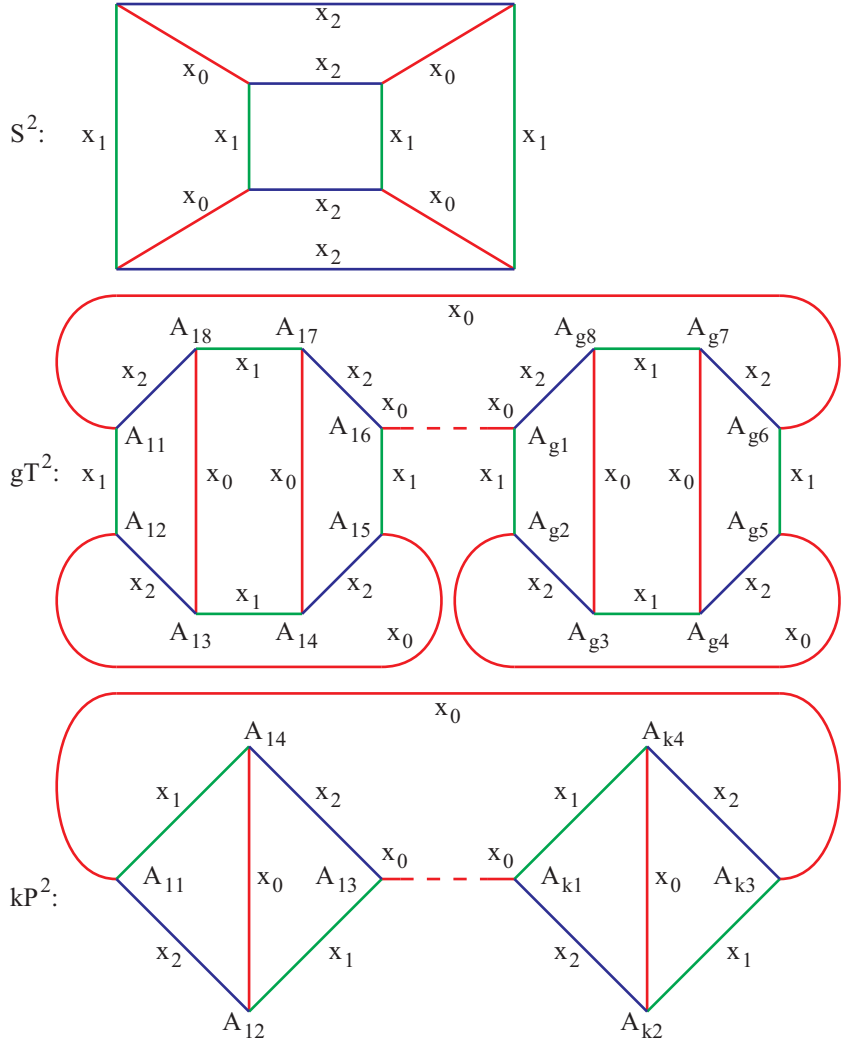}
 \caption{Graphs for coloring closed surfaces}\label{fg:surface}
\end{figure}

\begin{rem}
The first sentence is in fact a corollary of Theorem
\ref{th:colorable}. Notice that the graph used in
Fig.~\ref{fg:surface} for a real projective plane $P^2$ is exactly
the same as Fig.~\ref{fg:coloring} (b).
\end{rem}

\begin{proof}
The graphs for building up closed surfaces are shown in
Fig.~\ref{fg:surface}, where $\{x_0,x_1,x_2\}$ is a chosen basis for
$(\mathbb Z_2)^3$. By Theorem \ref{th:surface}, these graphs all
have $2$-skeletal expansions, and the corresponding topological
space will automatically become closed surfaces. Now let us verify
the topological type for these surfaces.

The graph for $S^2$ is just a simple exercise. Let us examine the
graph for $gT^2$ (the connected sum of $g$ copies of $T^2$). The
graph has $8g$ vertices, $12g$ edges and $2g+2$ bi-colored circles,
which are listed below:
$$\begin{tabular}{c|c|l}
 color & notation & \multicolumn{1}{c}{$2$-nest}\\
 \hline
 $\mathrm{Span}(x_0,x_1)$ & $\beta$ & $(A_{11}A_{12}A_{15}A_{16}\cdots
  A_{g1}A_{g2}A_{g5}A_{g6})$\\
 & $\gamma_1$ & $(A_{13}A_{14}A_{17}A_{18})$\\
 & \vdots & \vdots\\
 & $\gamma_g$ & $(A_{g3}A_{g4}A_{g7}A_{g8})$\\
 \hline
 $\mathrm{Span}(x_0,x_2)$ & $\xi$ &
  $(A_{11}A_{18}A_{13}A_{12}A_{15}A_{14}A_{17}A_{16}\cdots$\\
 & & \qquad$\cdots A_{g1}A_{g8}A_{g3}A_{g2}A_{g5}A_{g4}A_{g7}A_{g6})$\\
 \hline
 $\mathrm{Span}(x_1,x_2)$ & $\eta_1$ &
  $(A_{11}A_{12}A_{13}A_{14}A_{15}A_{16}A_{17}A_{18})$\\
 & \vdots & \vdots\\
 & $\eta_g$ & $(A_{g1}A_{g2}A_{g3}A_{g4}A_{g5}A_{g6}A_{g7}A_{g8})$\\
\end{tabular}$$

Suppose $(M,K)$ is its $2$-skeletal expansion. Then one can choose
an orientation on each of these bi-colored circles, such that when
two bi-colored circles share a common edge, they induce opposite
orientations on that edge. This implies that $M$ is an orientable
surface. Moreover,  its Euler characteristic is
$8g-12g+(2g+2)=2-2g$, so $M$ is indeed an orientable surface with
genus $g$.

As for the graph for $kP^2$ (i.e., a connected sum of $k$ copies of
$P^2$), it has $4k$ vertices, $6k$ edges and $k+2$ bi-colored
circles, which are listed below:
$$\begin{tabular}{c|c|l}
 color & notation & \multicolumn{1}{c}{$2$-nest}\\
 \hline
 $\mathrm{Span}(x_0,x_1)$ & $\beta$ & $(A_{11}A_{14}A_{12}A_{13}\cdots
  A_{k1}A_{k4}A_{k2}A_{k3})$\\
 \hline
 $\mathrm{Span}(x_0,x_2)$ & $\xi$ & $(A_{11}A_{12}A_{14}A_{13}\cdots
  A_{k1}A_{k2}A_{k4}A_{k3})$\\
 \hline
 $\mathrm{Span}(x_1,x_2)$ & $\eta_1$ & $(A_{11}A_{12}A_{13}A_{14})$\\
 & \vdots & \vdots\\
 & $\eta_k$ & $(A_{k1}A_{k2}A_{k3}A_{k4})$\\
\end{tabular}$$

Suppose $(M,K)$ is its $2$-skeletal expansion. The $2$-cells bounded
by $\beta$ and $\xi$ can not have opposite orientations on all their
common boundary arcs, so $M$ must be a non-orientable surface, while
its Euler characteristic is equal to $4k-6k+(k+2)=2-k$. Thus $M$ is
indeed a non-orientable surface with genus $k$.
\end{proof}

\begin{rem}
The graphs for $gT^2$ and $kP^2$ are actually constructed from the
graphs for $1T^2$ (torus) and $1P^2$ (real projective plane) by
doing \emph{connected sums}. Suppose $\Gamma_i$ is a regular
$3$-valent graph with ``pure'' coloring
$\alpha_i:E_{\Gamma_i}\to\{x_0,x_1,x_2\}$ and $(M_i,K_i)$ is its
$2$-skeletal expansion ($i=1$ or $2$). If $e_1,e_2$ are two edges in
$\Gamma_1,\Gamma_2$ respectively with the same color, then one can
cut them open, reconnect the four ends in the other way and form a
new regular $3$-valent graph $\Gamma$ with ``pure'' coloring
$\alpha:E_\Gamma\to\{x_0,x_1,x_2\}$. Suppose $(M,K)$ is its
$2$-skeletal expansion, then it is not difficult to see that $M$ is
exactly the connected sum of $M_1$ and $M_2$.
\end{rem}

\section{$(\mathbb Z_2)^4$-colorable closed $3$-manifolds}

Now let us consider the case $n=3$.

\begin{thm}\label{th:3-manifold}
Suppose $G=(\mathbb Z_2)^4$, $\Gamma$ is a connected regular
$4$-valent graph and $\alpha:E_\Gamma\to\mathrm{Hom}(G,\mathbb Z_2)$
is a ``good'' $G$-coloring. Then $(\Gamma,\alpha)$ has a
$3$-skeletal expansion $(M,K)$ (by Theorem \ref{th:manifold} $M$
must be a closed $3$-manifold) if and only if:
$$\#\{3\text{-nests}\}=\#\{2\text{-nests}\}-\#\{\text{vertices}\}.$$
\end{thm}

\begin{proof}
First let us use $\nu_k$ to denote (only in this proof) the number
of $k$-nests in $\Gamma$. Now suppose $(\Gamma,\alpha)$ has a
$3$-skeletal expansion $(M,K)$. Then by Theorem \ref{th:manifold},
$M$ is a closed $3$-manifold. So it is not difficult to see that
$\nu_0-\nu_1+\nu_2-\nu_3$ is equal to the Euler characteristic of
$M$, which must be $0$. Since $\Gamma$ is a regular $4$-valent
graph, $\nu_1=2\nu_0$ so we get $\nu_3=\nu_2-\nu_0$.

On the other hand, suppose that $\nu_3=\nu_2-\nu_0$. First choose a
$2$-skeletal expansion $(M_2,K_2)$ for $(\Gamma,\alpha)$ according
to Lemma \ref{lm:induction}. For each three edges $e_0,e_1,e_2$
sharing a common vertex, suppose $\Delta$ is the unique $3$-nest
containing them and suppose $F=F^2(e_0,e_1,e_2)$. Let
$\nu_k(\Delta)$ denote the number of $k$-nests in $\Delta$. Then the
Euler characteristic of $F$ is equal to
$\chi(\Delta)=\nu_0(\Delta)-\nu_1(\Delta)+\nu_2(\Delta)$, since $F$
is a $2$-skeletal expansion of $\Delta$. However, $F$ is in fact a
closed surface, thus $\chi(\Delta)\leq 2$. Because $\Delta$ is a
regular $3$-valent graph without loop edge,
$\nu_1(\Delta)=3\nu_0(\Delta)/2$. Thus,
$$\nu_2(\Delta)-\nu_0(\Delta)/2\leq 2.$$
Because any $k$ edges ($k=0,1,2$, or $3$) sharing a common vertex
determine a unique $k$-nest, we know that every vertex is the face
of exactly four $3$-nests, and every $2$-nest is the face of exactly
two $3$-nests. Therefore summing up the above inequality over all
$3$-nests, one obtains
$$2\nu_2-2\nu_0\leq 2\nu_3.$$

Our assumption implies that the inequalities here must all be
equalities. Consequently those $F^2(e_0,e_1,e_2)$ must all be
spheres. Hence by Lemma \ref{lm:induction}, one can fill these
spheres with some $3$-cells and get a $3$-skeletal expansion
$(M,K)$. By Theorem \ref{th:manifold}, $M$ must be a closed
$3$-manifold.
\end{proof}

\begin{figure}[ht]
 \includegraphics{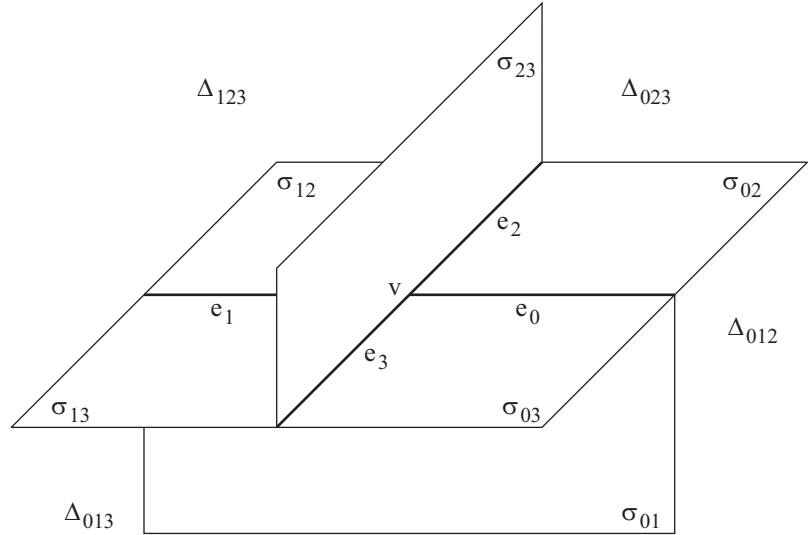}
 \caption{Local picture near a vertex of $3$-skeletal
  expansion}\label{fg:local}
\end{figure}

To understand intuitively why $M$ must be a closed $3$-manifold
(without invoking Theorem \ref{th:manifold}), we can consider the
local picture near a vertex $v$ (see Fig.~\ref{fg:local}). Choose
some local coordinates such that $v$ becomes $(0,0,0)$ and the
adjacent edges become (locally)
$e_0\rightsquigarrow\{(t,0,0)~|~0<t<1\}$,
$e_1\rightsquigarrow\{(t,0,0)~|-1<t<0\}$,
$e_2\rightsquigarrow\{(0,t,0)~|~0<t<1\}$,
$e_3\rightsquigarrow\{(0,t,0)~|-1<t<0\}$ respectively. Suppose
$\alpha(e_i)=x_i$.

For each $i\neq j\in\{0,1,2,3\}$, there is exactly one $2$-cell $\sigma_{ij}$
adjacent to both $e_i$ and $e_j$, let us glue it to $e_i$ and $e_j$ (locally)
as follows:
$$\begin{array}{l}
 \sigma_{01}\rightsquigarrow\{(t,0,s)~|-1<t<1,-1<s<0\}\\
 \sigma_{02}\rightsquigarrow\{(t,s,0)~|~0<t<1,0<s<1\}\\
 \sigma_{03}\rightsquigarrow\{(t,s,0)~|~0<t<1,-1<s<0\}\\
 \sigma_{12}\rightsquigarrow\{(t,s,0)~|-1<t<0,0<s<1\}\\
 \sigma_{13}\rightsquigarrow\{(t,s,0)~|-1<t<0,-1<s<0\}\\
 \sigma_{23}\rightsquigarrow\{(0,t,s)~|-1<t<1,0<s<1\}\\
\end{array}$$

For each $i\neq j\neq k\in\{0,1,2,3\}$, there is exactly one $3$-cell
$\Delta_{ijk}$ adjacent to three edges $e_i$, $e_j$ and $e_k$, and we can glue
it to $e_i$, $e_j$ and $e_k$ (locally) as follows:
$$\begin{array}{l}
 \Delta_{012}\rightsquigarrow\{(t,s,u)~|-1<t<1,0<s<1,-1<u<0\}\\
 \Delta_{013}\rightsquigarrow\{(t,s,u)~|-1<t<1,-1<s<0,-1<u<0\}\\
 \Delta_{023}\rightsquigarrow\{(t,s,u)~|~0<t<1,-1<s<1,0<u<1\}\\
 \Delta_{123}\rightsquigarrow\{(t,s,u)~|-1<t<0,-1<s<1,0<u<1\}\\
\end{array}$$

The possibility of this arrangement implies that every vertex $v$
has a neighborhood in $M$ homeomorphic to $\mathbb R^3$. Thus the
combinatorial nature of $(M,K)$ implies that $M$ is a closed
$3$-manifold.

On the other hand, one can also use this canonical method to verify
the $G$-colorability of a general $3$-manifold, by finding cell
decomposition that locally looks like the above canonical picture
near each vertex.

\begin{thm}
Every closed $3$-manifold is $(\mathbb Z_2)^4$-colorable. Moreover, one can
use only four $3$-cells to do this.
\end{thm}

\begin{rem}
Notice that the first half of the statement is also a corollary of
Theorem \ref{th:colorable}, since every closed $3$-manifold has a
(in fact unique up to diffeomorphism) differential structure (see
e.g. \cite{M}), which implies that every closed $3$-manifold is a
closed combinatorial manifold. Notice also that for $3$-manifolds
the graph may \emph{not} be planar anymore, and may contain
double-edges.
\end{rem}

\begin{proof}
Recall the idea of Heegaard splitting (see e.g. \cite{H}): suppose
$H^\pm$ are two genus-$g$ handle-bodies (here we require $g>0$),
$F^\pm=\partial H^\pm$ are two genus-$g$ closed orientable surfaces,
and $f:F^-\to F^+$ is a homeomorphism, then one can glue $H^-$ to
$H^+$ along their boundaries using $f$ and get a closed $3$-manifold
$M=H^+\cup_f H^-$. Every closed $3$-manifold $M$ can be obtained in
this way and the corresponding $(H^+,H^-,f)$ is called a
\emph{Heegaard splitting} for $M$.

\begin{figure}[ht]
 \includegraphics{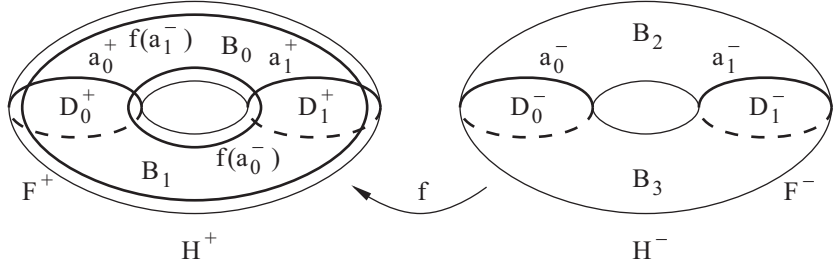}
 \caption{A Heegaard splitting for $S^3$}\label{fg:heegaard}
\end{figure}

Each genus-$g$ handle-body $H$ can be further decomposed into two balls
$B_\pm(H)$ in a canonical way along $(g+1)$ disjoint non-separating discs
$D_0(H),\ldots,D_g(H)$. Let us decompose $H^+$ and $H^-$ in this way. For
notational convenience, let $D_i^\pm=D_i(H^\pm)$
$$B_0=B_+(H^+),B_1=B_-(H^+),B_2=B_+(H^-),B_3=B_-(H^-),$$
and let $a^\pm_i=\partial(D_i^\pm)\subseteq\partial H^\pm$ (see
Fig.~\ref{fg:heegaard}).

Notice that the topological type of $M=H^+\cup_f H^-$ remains
unchanged when changing $f$ isotopically, while these $f(a^-_j)$
will be moved around on $F^+$ when doing so. Clearly one can
isotopically adjust $f$ a little bit, such that each $f(a^-_j)$
intersects these $a^+_i$ transversely in a finite number of points.
Moreover, one can isotopically adjust $f$ a little bit further,
adding some redundant turns to those $f(a^-_j)$, such that for any
$i,j\in\{0,1,\ldots,g\}$,
$$a^+_i\cap(f(a^-_0)\cup\cdots\cup f(a^-_g))\neq\emptyset,
 ~(a^+_0\cup\cdots\cup a^+_g)\cap f(a^-_j)\neq\emptyset,$$
and all these $a^+_i$ and $f(a^-_j)$ together cut $F^+$ into small
\emph{disc} regions.

Now the points in $(a^+_0\cup\cdots\cup
a^+_g)\cap(f(a^-_0)\cup\cdots\cup f(a^-_g))$ cut these $a^+_i$ and
$f(a^-_j)$ into edges, and these edges together form a finite
connected regular $4$-valent graph $\Gamma\subseteq F^+$. Suppose
$\{t_0, t_1, t_2,t_3\}$ is a basis for $G$, and $\{x_0, x_1,
x_2,x_3\}$ is the dual basis in $\mathrm{Hom}(G,\mathbb Z_2)$,
namely $x_i(t_j)=\delta_{ij}$. Define a characteristic function
$\lambda:B_i\mapsto t_i$ and then define the edge coloring
$\alpha:E_\Gamma\to\mathrm{Hom}(G,\mathbb Z_2)$ as follows:
$\alpha(e)(\lambda(B_i))=0$ if and only if $e\subset\partial B_i$.
In fact, each $e\in E_\Gamma$ lies within the boundary of exactly
three of $B_0, B_1, B_2,B_3$. If the fourth ball is $B_i$, then
$\alpha(e)=x_i$. It is easy to see that the set of colors of the
four edges around each vertex of $\Gamma$ is exactly $\{x_0, x_1,
x_2,x_3\}$, so $\alpha$ is a ``pure'' coloring.

Let us calculate the $3$-skeletal expansion of $(\Gamma,\alpha)$.
Since $\alpha$ is a ``pure'' coloring, the set of all 2-nests of
$(\Gamma,\alpha)$ consists of exactly all the bi-colored circles of
$\Gamma$. By the construction of $\Gamma$ as described above, it is
easy to see that each bi-colored circle of $\Gamma$ is either the
boundary of a disc region of
$$F^+\setminus(a^+_0\cup\cdots\cup a^+_g\cup f(a^-_0)\cup\cdots\cup
 f(a^-_g))$$
or one of those $a^+_i$ (bounds $D_i^+$ in $H^+$) or $f(a^-_j)$
(bounds $D_j^-$ in $H^-$). Let those $D_i^\pm$ and disc regions be
the $2$-cells. Then one may construct a $2$-skeletal expansion
$(M_2,K_2)$ of $(\Gamma, \alpha)$, where $M_2$ is obtained by gluing
each $D_i^-$ to $F^+\cup D_0^+\cup\cdots\cup D_g^+$ along the
boundary, such that $a_i^-$ is identified with $f(a_i^-)$. Since an
edge $e\subseteq\partial B_i$ if and only if $\alpha(e)\neq x_i$,
all edges in each $\partial B_i$ form a connected regular $3$-valent
subgraph $\Delta_i\subseteq\Gamma$ with $\dim_\alpha\Delta_i=3$, and
$\Delta_0,\Delta_1,\Delta_2,\Delta_3$ are all the $3$-nests. Hence
by adding $B_0,B_1,B_2,B_3$ as the $3$-cells, one obtains a
$3$-skeletal expansion $(M_3,K_3)$. Clearly $M_3=\bigcup_{\sigma\in
K_3}\sigma$ is just the original closed $3$-manifold $M$.
\end{proof}

\begin{figure}[ht]
 \includegraphics{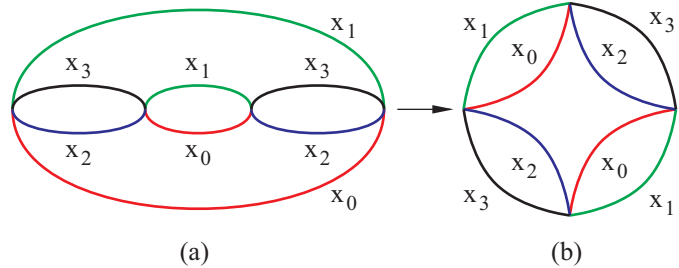}
 \caption{Colored graph for $S^3$}\label{fg:graph}
\end{figure}

As an example, Fig. \ref{fg:graph} (a) shows the colored graph corresponding
to the Heegaard splitting of $S^3$ in Fig. \ref{fg:heegaard}. For a more fancy
picture of this graph see Fig. \ref{fg:graph} (b).

\begin{figure}[ht]
 \includegraphics{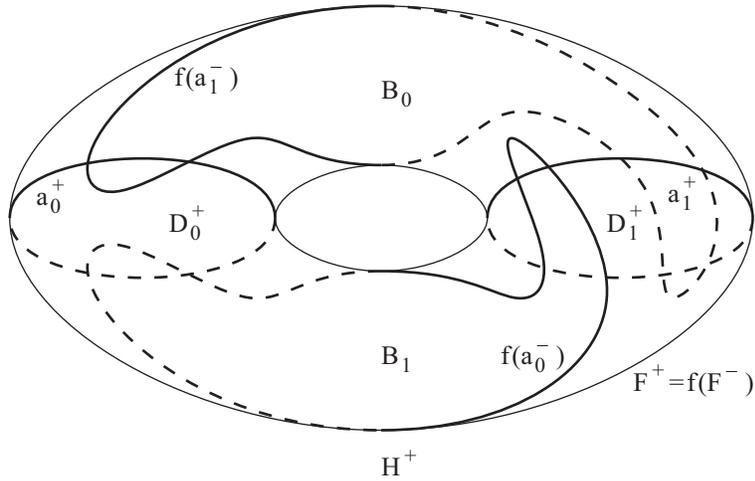}
 \caption{A Heegaard splitting for $S^2\times S^1$}\label{fg:example}
\end{figure}

\begin{figure}[ht]
 \includegraphics{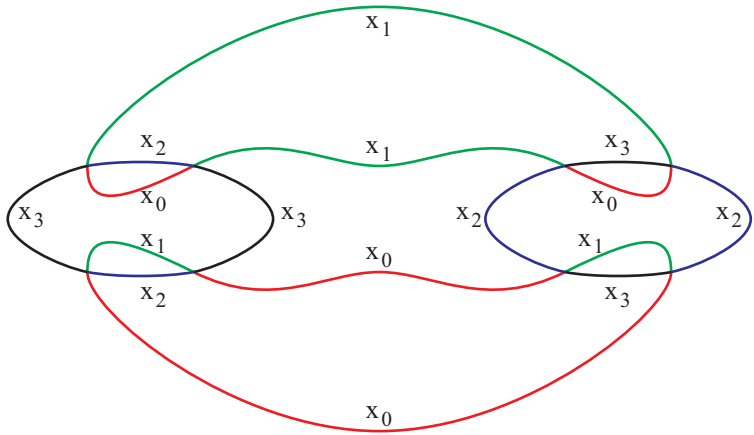}
 \caption{Colored graph for $S^2\times S^1$}\label{fg:examplegraph}
\end{figure}

Another example is $S^2\times S^1$. Fig. \ref{fg:example} shows a nice
Heegaard splitting, such that every small region on $F^+$ is a disc (the
picture for $H^-$ is omitted), and Fig. \ref{fg:examplegraph} is the
corresponding colored graph.

\begin{rem}
Similarly to the case of simple convex polytopes, the colorings here
are also induced by characteristic functions on the $n$-dimensional
faces. However, the cell decompositions here are not dual to any
triangulation, so this is different from the case of simple convex
polytopes or the graphs constructed in the proof of Theorem
\ref{th:colorable}. In fact a triangulation for a closed
$3$-manifold has at least two $3$-simplices, so it has at least five
vertices, which implies that its dual regular cell decomposition has
at least five $3$-cells.
\end{rem}

\section{High dimensional cases}

\begin{proof}[Proof of Theorem \ref{th:manifold}]
Because of the combinatorial nature of $(M,K)$, for any cell
$\sigma\in K$ with dimension $\geq 1$ and any two points $u,v$
inside $\sigma$ (away from $\partial\sigma$), $u$ has a neighborhood
homeomorphic to $\mathbb R^n$ if and only if $v$ does so. Therefore
we only need to verify that every point in $M$ sufficiently close to
a vertex has a neighborhood homeomorphic to $\mathbb R^n$. Now
choose a vertex $v$, we will show that $v$ has a neighborhood
homeomorphic to the space
$$C_n=\{(t_0,\ldots,t_n)\in\mathbb R^{n+1}|~t_0\cdots t_n=0,0\leq
 t_0<1,\ldots,0\leq t_n<1\},$$
which is homeomorphic to $\mathbb R^n$.

The argument is quite similar to the case $n=3$ in section 3. Choose
a local coordinate such that $v$ becomes $(0,\ldots,0)$ and the
adjacent edges become (locally)
$e_0\rightsquigarrow\{(t,0,\ldots,0)~|~0<t<1\},
e_1\rightsquigarrow\{(0,t,0,\ldots,0)~|~0<t<1\},\ldots,
e_n\rightsquigarrow\{(0,\ldots,0,t)~|~0<t<1\}$ respectively. Suppose
$\alpha(e_i)=x_i$.

For each $i_1<\cdots<i_k\in\{0,\ldots,n\}$, there is a unique $k$-nest
$\Delta_{i_1\cdots i_k}$ containing $e_{i_1},\ldots,e_{i_k}$, and by Lemma
\ref{lm:regular} these $\Delta_{i_1\cdots i_k}$ enumerate all $k$-nests
containing $v$. Suppose $\sigma_{i_1\cdots i_k}$ is the $k$-cell corresponding
to $\Delta_{i_1\cdots i_k}$. By an induction on $k$ one can adjust the local
coordinates such that each of these $\sigma_{i_1\cdots i_k}$ becomes (locally)
$$\{(t_0,\dots,t_n)~|~0<t_{i_1}<1,\ldots,0<t_{i_k}<1,t_j=0\text{ for all
 other }j\}$$

The possibility of this arrangement implies that every vertex $v$
has a neighborhood in $M$ homeomorphic to the cone $C_n$, which is
homeomorphic to $\mathbb R^n$. Thus the combinatorial nature of
$(M,K)$ implies that every point in $M$ has a neighborhood
homeomorphic to $\mathbb R^n$. $M$ is connected since $\Gamma$ is
connected, and $M$ is compact because it has only finitely many
cells. Therefore $M$ is a closed $n$-manifold.
\end{proof}

Before proving Theorem \ref{th:colorable}, let us recall some
notation. If a cell $\sigma$ is a face of cell $\sigma'$ and
$\sigma\neq\sigma'$, we write $\sigma\prec\sigma'$. Now suppose that
$K$ is a finite regular cell decomposition of $M$, then the
\emph{barycentric subdivision} $\mathrm{Sd}~K$ is a simplicial
complex, such that every $k$-simplex in $\mathrm{Sd}~K$ corresponds
to a sequence of cells $\sigma_0\prec\cdots\prec\sigma_k\in K$, we
denote such a simplex by $[\sigma_0\cdots\sigma_k]$, and call it a
\emph{flag}. Notice that such a simplicial complex always exists and
is a triangulation of $M$. When $K$ itself is a triangulation of
$M$, $[\sigma_0\cdots\sigma_k]$ is just the convex hull spanned by
$\tilde\sigma_0,\ldots,\tilde\sigma_k$, where $\tilde\sigma$ denotes
the barycenter of $\sigma$.

For each cell $\sigma\in K$, the \emph{link complex} is
$$\mathrm{Lk}~\sigma=\{[\sigma_0\cdots\sigma_k]\in\mathrm{Sd}~K~|~
 \sigma\prec\sigma_0\prec\cdots\prec\sigma_k\}.$$
If $M$ is connected and has a finite regular cell decomposition $K$, and
$$\forall\sigma\in K,~|\mathrm{Lk}~\sigma|\cong S^{n-\dim\sigma-1}$$
($|\mathrm{Lk}~\sigma|$ is the union of all cells in $\mathrm{Lk}~\sigma$),
then $M$ is called an \emph{$n$-dimensional closed combinatorial manifold}.
Notice that if $K$ satisfies this condition, then $\mathrm{Sd}~K$ also
satisfies this condition.

Finally for each cell $\sigma\in K$, the \emph{dual block} is
$$\mathfrak D\sigma=\{[\sigma_0\cdots\sigma_k]\in\mathrm{Sd}~K~|~
 \sigma=\sigma_0\prec\cdots\prec\sigma_k\}.$$
If $K$ satisfies the previous condition, then these dual blocks are
all cells (in fact $\overline{\mathfrak D\sigma}$ is the topological
cone on $\mathrm{Lk}~\sigma$), and they also form a finite regular
cell decomposition of $M$, we denote it by $\mathfrak DK$. In
particular, the $1$-skeleton of this dual cell decomposition is a
regular graph. It is clear that $\mathfrak D:K\to\mathfrak DK$ is a
one-to-one correspondence \emph{reversing face relations}, and
$\dim\mathfrak D\sigma=n-\dim\sigma$.

\begin{proof}[Proof of Theorem \ref{th:colorable}]
Suppose that $K'$ is a finite regular cell decomposition of $M$
satisfying the condition
$$\forall\sigma\in K',~|\mathrm{Lk}~\sigma|\cong
 S^{n-\dim\sigma-1}.$$
Let $K''=\mathrm{Sd}~K'$ be its barycentric subdivision. Then $K''$
also satisfies the above condition, thus one can define a dual
regular cell decomposition $K=\mathfrak DK''$. Now let us define a
coloring $\alpha$ for the $1$-skeleton of $K$ (denoted by $\Gamma$)
as follows.

Let $x_0,\ldots,x_n$ be a linear basis for $\mathrm{Hom}(G,\mathbb Z_2)$. For
each $1$-cell $\sigma=\mathfrak D_{K''}[\sigma_0\cdots\sigma_{n-1}]\in\Gamma$
(here $\sigma_0\prec\cdots\prec\sigma_{n-1}\in K'$), suppose that
$$\{0,\ldots,n\}\setminus\{\dim\sigma_0,\ldots,\dim\sigma_{n-1}\}=\{k\},$$
then let $\alpha(\sigma)=x_k$. Because $\mathfrak D_{K''}$ reverses
face relations, for each vertex $v\in K$, say $v=\mathfrak
D_{K''}[\sigma_0\cdots\sigma_n]$ (here
$\sigma_0\prec\cdots\prec\sigma_n\in K'$), there are exactly $(n+1)$
edges in $K$ adjacent to $v$, namely
$$e_k=\mathfrak
 D_{K''}[\sigma_0\cdots\sigma_{k-1}\sigma_{k+1}\cdots\sigma_n],~k=0,\ldots,n.$$
Since $\alpha(e_k)=x_k$, this implies that $(\Gamma,\alpha)$ is a finite
regular $(n+1)$-valent graph with ``pure'' $G$-coloring (see Fig
\ref{fg:dual}).

\begin{figure}[ht]
 \includegraphics{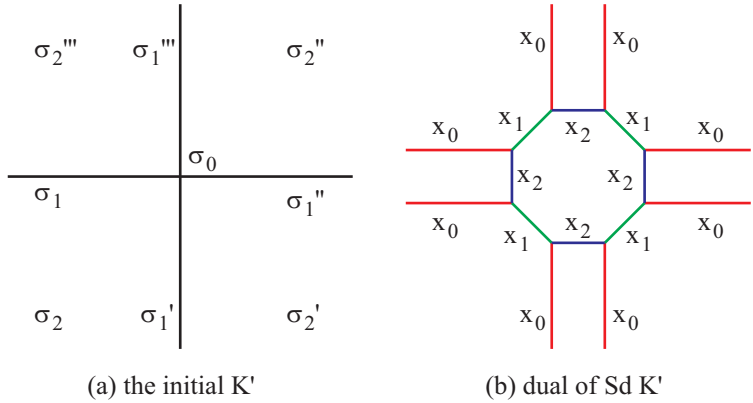}
 \caption{The dual graph of barycentric subdivision}\label{fg:dual}
\end{figure}

Now for any $m$-cell $\Delta=\mathfrak
D_{K''}[\sigma_{j_0}\cdots\sigma_{j_{n-m}}]$ (here
$\sigma_{j_0}\prec\cdots\prec\sigma_{j_{n-m}}\in K'$ and
$\dim\sigma_{j_k}=j_k$), let $\kappa(\Delta)$ be the $1$-skeleton of
$\partial\Delta$. Then a vertex $v'=\mathfrak
D_{K''}[\sigma'_0\cdots\sigma'_n]$ (here
$\sigma'_0\prec\cdots\prec\sigma'_n\in K'$) is a vertex of $\kappa(\Delta)$ if
and only if
$$\sigma'_{j_0}=\sigma_{j_0},\ldots,\sigma'_{j_{n-m}}=\sigma_{j_{n-m}},$$
and an edge $e'=\mathfrak
D_{K''}[\sigma'_0\cdots\sigma'_{k-1}\sigma'_{k+1}\cdots\sigma'_n]$ adjacent to
$v'$ is an edge of $\kappa(\Delta)$ if and only if
$k\not\in\{j_0,\ldots,j_{n-m}\}$, namely
$\alpha(e')\in\{x_{i_1},\ldots,x_{i_m}\}$
$$(i_1<\cdots<i_m\in\{0,\ldots,n\}\setminus\{j_0,\ldots,j_{n-m}\}).$$
This implies that $\kappa(\Delta)$ is a connected regular $m$-valent subgraph
of $\Gamma$ with edge-colors in $\{x_{i_1},\ldots,x_{i_m}\}$, hence
$\kappa(\Delta)$ is an $m$-nest in $K(\Gamma,\alpha)$.

Therefore the correspondence $\kappa:K\to K(\Gamma,\alpha)$, sending
each cell in $K$ to the $1$-skeleton of its boundary complex, is a
well defined map which preserves dimensions and face relations.
$\kappa$ is also a one-to-one correspondence. In fact, for any
$m$-nest $\Delta$, suppose $v=\mathfrak
D_{K''}[\sigma_0\cdots\sigma_n]$ is a vertex in $\Delta$ and
$e_{i_1},\ldots,e_{i_m}$ are the $m$ edges in $\Delta$ adjacent to
$v$, with each $e_{i_k}=\mathfrak
D_{K''}[\sigma_0\cdots\sigma_{i_k-1}\sigma_{i_k+1}\cdots\sigma_n]$.
Then any inverse image of $\Delta$ under $\kappa$ must contain $v$
and these edges, and there does exists a unique $m$-cell containing
$e_{i_1},\ldots,e_{i_m}$, namely $\mathfrak
D_{K''}[\sigma_{j_0}\cdots\sigma_{j_{n-m}}]$ (here
$j_0<\cdots<j_{n-m}\in\{0,\ldots,n\}\setminus\{i_1,\ldots,i_m\})$.
Clearly $\mathfrak D_{K''}[\sigma_{j_0}\cdots\sigma_{j_{n-m}}]$ is
the one and only cell that is sent to $\Delta$ by $\kappa$. Thus
$\kappa$ must be a one-to-one correspondence.
\end{proof}

As an example of calculation, let us reproduce the ``cubic'' graph for $S^n$
as a $\mathfrak D(\mathrm{Sd}~K')$. Take a regular cell decomposition $K'$ for
$S^n$ with two $i$-cells at each dimension $i\in\{0,1,\ldots,n\}$:
$$\begin{array}{l}
 \sigma_i^{(1)}=\{(t_0,\ldots,t_n)\in\mathbb R^{n+1}|\\
 \qquad t_0^2+\cdots+t_i^2=1,t_i>0,t_{i+1}=\cdots=t_n=0\};\\
 \sigma_i^{(2)}=\{(t_0,\ldots,t_n)\in\mathbb R^{n+1}|\\
 \qquad t_0^2+\cdots+t_i^2=1,t_i<0,t_{i+1}=\cdots=t_n=0\}.
\end{array}$$
Notice that $\sigma_i^{(\epsilon_i)}\prec\sigma_j^{(\epsilon_j)}$ if
and only if $i<j$. Therefore $K=\mathfrak D(\mathrm{Sd}~K')$
contains $2^{n+1}$ vertices, each of which has the form
$$\mathfrak D[\sigma_0^{(\epsilon_0)}\cdots\sigma_n^{(\epsilon_n)}],~
 \epsilon_0,\ldots,\epsilon_n\in\{1,2\}.$$
To give the graph a more fancy appearance, identify each such vertex
with
$$v(\epsilon)=((-1)^{\epsilon_0},\ldots,(-1)^{\epsilon_n})\in \mathbb
 R^{n+1}.$$
Two vertices $v(\epsilon)$ and $v(\tau)$ are connected by an edge in $K$ of
the form
$$\mathfrak D[\sigma_0^{(\epsilon_0)}\cdots\sigma_{i-1}^{(\epsilon_{i-1})}
 \sigma_{i+1}^{(\epsilon_{i+1})}\cdots\sigma_n^{(\epsilon_n)}]$$
(colored by $x_i$) if and only if $\tau_i\neq\epsilon_i$ while for
all other $j$, $\tau_j=\epsilon_j$. This implies that we can also
identify each such edge with the line segment
$$\{(t_0,\ldots,t_n)\in \mathbb R^{n+1}|
 -1\leq t_i\leq 1,~\forall j\neq i,t_j=(-1)^{\epsilon_j}\}.$$
Hence one obtains the $1$-skeleton of the $(n+1)$-dimensional cube
$$\mathcal C=[-1,1]^{n+1}\subset \mathbb R^{n+1},$$
and an edge of this graph is colored by $x_i$ if and only if it is
parallel to the $t_i$-axis. When $n=2$, this equals to the graph for
$S^2$ in Fig. \ref{fg:surface}.

\begin{rem}
There is in fact a more general way of doing skeletal expansions,
which we call the ``generalized skeletal expansion''. For example,
take the colored graph in Fig. \ref{fg:coloring} (a), one can also
glue M\"obius bands instead of discs to each of the bi-colored
circles, and then the resulting manifold will be a closed
non-orientable surface with genus $4$.

To describe this expansion more accurately, we will use an induction
on the dimension. The generalized $1$-skeletal expansion of a
regular graph $\Gamma$ with ``good'' coloring $\alpha$ is just
$\Gamma$ itself. For $k>1$, a generalized $k$-skeletal expansion of
$(\Gamma,\alpha)$ is a space $M^k$ together with a defining sequence
$\mathfrak E^k$ of gluing operations such that $\Gamma$ can be
extended to $M^k$ by $\mathfrak E^k$ and each $\ell$-nest with
$\ell\leq k+1$ can be extended to a closed $(\ell-1)$-dimensional
manifold. If $(\Gamma,\alpha)$ admits a generalized $k$-skeletal
expansion $(M^k,\mathfrak E^k)$, then $\mathfrak E^k$ also induces a
generalized $k$-skeletal expansion for each $(k+1)$-nest of
$(\Gamma,\alpha)$. By $(F^k,\mathfrak E^k)$ we denote the disjoint
union of the generalized $k$-skeletal expansions of all
$(k+1)$-nests in $\Gamma$. Suppose that $(\Gamma,\alpha)$ admits a
generalized $k$-skeletal expansion $(M^k,\mathfrak E^k)$ such that
$F^k$ is the boundary of a $(k+1)$-dimensional compact manifold
$\tilde M^{k+1}$ (may not be connected). Then one can glue $\tilde
M^{k+1}$ onto $M^k$ to form a new space $M^{k+1}$. Add this
operation to $\mathfrak E^k$ to form a new sequence of operations
$\mathfrak E^{k+1}$, and we call $(M^{k+1},\mathfrak E^{k+1})$ a
\emph{generalized $(k+1)$-skeletal expansion} of $(\Gamma,\alpha)$.
Clearly, in each dimension, ordinary skeletal expansions (if they
exist) are uniquely determined by the colored graphs up to cell
isomorphisms, while the generalized skeletal expansions are not so.
\end{rem}

\begin{figure}[ht]
 \includegraphics{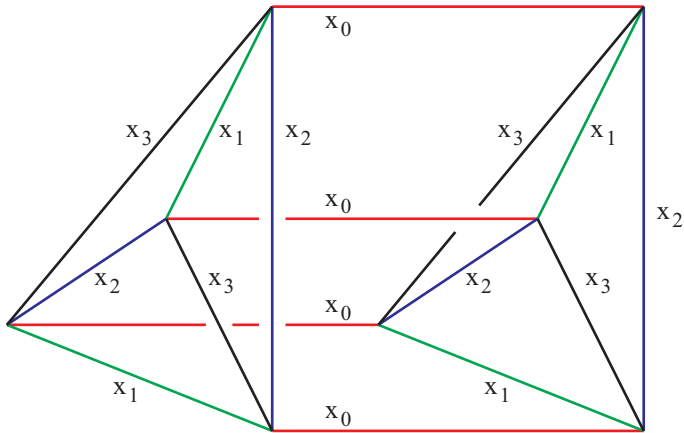}
 \caption{Graph admitting generalized $3$-skeletal
  expansion}\label{fg:general}
\end{figure}

As an example, let us study the colored regular $4$-valent graph
$\Gamma$ shown in Fig. \ref{fg:general}. It admits a $2$-skeletal
expansion but admits no $3$-skeletal expansion, since by Theorem
\ref{th:3-manifold}
$$\#\{3\text{-nests}\}-\#\{2\text{-nests}\}+\#\{\text{vertices}\}=5-12+8\neq
 0.$$
In fact $\Gamma$ has five $3$-nests. When doing (ordinary)
$2$-skeletal expansion on their disjoint union, two of them (both
colored by $x_1 x_2 x_3$) yield $P^2$, while the other three
(colored by $x_0 x_1 x_2$, $x_0 x_1 x_3$ and $x_0 x_2 x_3$
respectively) yield $S^2$. However, one can still glue three
$3$-cells to those $S^2$, glue a $P^2\times I$ to the two $P^2$, and
finally glue the boundary surfaces of all these blocks together
according to the chosen (ordinary) $2$-skeletal expansion. The
resulting space of this generalized $3$-skeletal expansion is
homeomorphic to
 $P^2\times S^1$. This graph $\Gamma$ even has
a generalized $4$-skeletal expansion, since $P^2\times S^1$ is the
boundary of the compact $4$-manifold  $P^2\times B^2$.

For the generalized skeletal expansion, we would like to propose the
following conjecture:

\begin{cnj}
Suppose that $\Gamma$ is a regular $(n+1)$-valent graph with a
``good'' $(\mathbb Z_2)^{n+1}$-coloring $\alpha$. Then
 $(\Gamma,\alpha)$ always admits a generalized $(n+1)$-skeletal
expansion $(M,\mathfrak E)$ such that $M$ is a compact
$(n+1)$-manifold with boundary.
\end{cnj}

\section{Reconstruction of manifolds with $(\mathbb
Z_2)^{n+1}$-actions}

Finally, let us consider the realization problem of reconstructing
$(n+1)$-dimensional $G$-manifolds from $(n+1)$-valent $G$-colored
graphs, where $G=(\mathbb Z_2)^{n+1}$.

In the previous section, we showed that under certain conditions, a
$G$-colored graph $(\Gamma,\alpha)$  admits a skeletal expansion
into a closed manifold $M$.  If $M$ can be used as the boundary of a
simple convex polytope $P$, then $\alpha$ determines a
characteristic function $\lambda$ on $P$. Furthermore, by the
reconstruction technique of small covers in \cite{DJ}, one can use
$\lambda$ and the product bundle $P\times G$ to construct a small
cover $X$ such that the moment graph of $X$ is exactly the given
$G$-colored graph $(\Gamma,\alpha)$. This observation provided much
insight into the study of the realization problem.

Now let us prove Theorem~\ref{realization}.

\begin{proof}[Proof of Theorem~\ref{realization}]
Let $\Gamma$ be a finite connected regular $(n+1)$-valent graph with
a ``good'' $G$-coloring $\alpha$. Suppose that $(\Gamma, \alpha)$
admits an $n$-skeletal expansion $(M,K)$. Then by
Theorem~\ref{th:manifold}, $M$ is an $n$-dimensional closed
manifold.

First, let us prove the special case in which $M$ is the boundary of
a compact $(n+1)$-dimensional manifold $X_0$. Suppose that
$\Delta\subset \Gamma$ is a $k$-nest with $k\leq n$. We recall that
the color of an edge $e$ is a nontrivial element $\alpha(e)$ in
$\mathrm{Hom}(G,\mathbb Z_2)$. Let $E_\Delta$ denote the set of all
edges in $\Delta$. Then $\mathfrak{g}_\Delta=\cap_{e\in
E_\Delta}\ker\alpha(e)$ is a co-rank $k$ subgroup of $G$. Notice
that $\mathfrak{g}_\Delta=G$ if $\Delta$ is a 0-nest (i.e., a single
vertex). On the other hand, the $(k-1)$-skeletal expansion of
$\Delta$ bounds a ball, and it corresponds to a $k$-dimensional open
cell $C_\Delta$ in $M$. Moreover, for each point $x\in M$, there is
a unique nest $\Delta$ such that $x\in C_\Delta$.

Now suppose that $X_0$ admits a trivial action of $G$. Consider the
trivial principal $G$-bundle $\xi=X_0\times G$ over $X_0$, where the
bundle projection is denoted by $\pi: \xi\longrightarrow X_0,
\pi(x,h)=x$, and the action of $G$ on $\xi$ is defined by $\phi:
G\longrightarrow \mathrm{End}(\xi)$ as follows:
$\phi(g)(x,h)=(x,gh)$. Then we define an equivalence relation $\sim$
on $\xi$ as follows: $p\sim p'$ if and only if
\begin{enumerate}
 \item[(1)]{$\pi(p)=\pi(p')\in X_0\backslash M$; or}
 \item[(2)] {there exists a nest $\Delta$ and some $g\in\mathfrak{g}_\Delta$ such that
 $\pi(p)=\pi(p')\in C_\Delta$ and $p=\phi(g)p'$.}
\end{enumerate}
By $[p]$ we denote the equivalence class of $p\in \xi$. Set
$X=\xi/\sim$.

Clearly $X$ is obtained by gluing $2^{n+1}$ copies of $X_0$ together
along their boundaries. To prove that $M$ is indeed an
$(n+1)$-dimensional manifold, it suffices to show that each point
$[p]$ with $\pi(p)\in M$ has a neighborhood homeomorphic to $\mathbb
R^{n+1}$. Let $[p]\in X$ with $\pi(p)\in M$. Then there is a
$k$-nest $\Delta$ such that $\pi(p)\in C_\Delta$. Since all points
inside a component of $\pi^{-1}(C_\Delta)$ can have homeomorphic
neighborhoods, we need only to prove that $[p]$ has a neighborhood
homeomorphic to $\mathbb R^{n+1}$ when $\pi(p)$ is very close to a
vertex.

Without loss of generality, one assumes that $v=\pi(p)$ is a vertex
of $\Gamma$. Then there are $n+1$ $n$-nests containing $v$, denoted
by $\Delta_0, ..., \Delta_n$ respectively. Choose a local coordinate
chart $(U, \varphi)$ near $v$ in $X_0$ such that $v$ is mapped into
the origin of $\mathbb R^{n+1}$, each $C_{\Delta_j}\cap U$ is mapped
into an $n$-dimensional open cone
$$\{(x_0,...,x_n)\in\mathbb R^{n+1}| x_j=0, x_i>0 \text{ when }
 i\not=j\},$$
and $U$ is mapped into an $(n+1)$-dimensional cone
$$\{(x_0,...,x_n)\in\mathbb R^{n+1}| x_i\geq 0, i=0,...,n \}.$$
Each $\mathfrak{g}_{\Delta_j}$ is a subgroup of rank 1, so we can
write $\mathfrak{g}_{\Delta_j}=\{\text{id}, g_j\}$. Obviously,
$\{g_0,...,g_n\}$ forms a basis of $G$. Then one can define a linear
action $\rho:G\longrightarrow \mathrm{End}(\mathbb R^{n+1})$ by
$$\rho(g_j)(x_0,...,x_j,...,x_n)=(x_0,...,-x_j,...,x_n).$$
Using this action, one can extend $\varphi$ into a continuous
surjection $\widetilde{\varphi}: \pi^{-1}(U)\longrightarrow\mathbb
R^{n+1}$ defined by
$$\widetilde{\varphi}(x, g)=\rho(g)(\varphi(x)).$$
It is not difficult to see that
$\widetilde{\varphi}(q)=\widetilde{\varphi}(q')$ if and only if
$q\sim q'$. Therefore, $\widetilde{\varphi}$ induces a homeomorphism
from a neighborhood of $[p]$ to $\mathbb R^{n+1}$.

The action $\phi$ of $G$ on $\xi$ also induces a natural action
$\widetilde{\phi}$ on $X$. If $\Delta$ is a $k$-nest, let
$\Phi_\Delta=\pi^{-1}(C_\Delta)/\sim$, then $G$ acts on
$\Phi_\Delta$ with kernel
$\ker\widetilde{\phi}=\mathfrak{g}_\Delta$. In fact, $\Phi_\Delta$
contains $2^k$ copies of $C_\Delta$ and the action of $G$ on
$\Phi_\Delta$ permutes these copies. Also, the fixed point set of
$\mathfrak{g}_\Delta$ acting on $X$ contains $\cup_{\Delta'\subseteq
\Delta}\Phi_{\Delta'}$ as a component. In particular, for an edge
$e$ with two ends $p, q$, the fixed point set of
$\mathfrak{g}_e=\ker\alpha(e)$ acting on $X$ contains a circle which
is the union of $\Phi_e$ (two arcs) and $\Phi_p, \Phi_q$ (two fixed
points). This also implies that the moment graph of the $G$-manifold
$X$ is just $(\Gamma, \alpha)$.

Now let us consider the general case in which $M$ is not the
boundary of a compact $(n+1)$-dimensional manifold. Choose an
$n$-nest $\Delta$, by the construction of $M$, then $\Delta$
corresponds to an open $n$-cell $C_\Delta$ of $(M, K)$ and the
$n$-skeletal expansion of $\Delta$ is exactly the closure
$\bar{C}_\Delta$ of $C_\Delta$, which is an $n$-disc. Now let us
remove a very small open $n$-ball $B$ in the interior of
$\bar{C}_\Delta$. Next, taking another copy $M'$ of $M$ and
forgetting those combinatorial structures $(\Gamma, \alpha)$ and $K$
on it, one randomly removes a small open $n$-ball $B'$ in $M'$. Then
one glues $M\backslash B$ and $M'\backslash B'$ together along their
boundaries. This is  actually  a connected sum between two copies of
$M$, so the resulting manifold $N$ becomes the boundary of an
$(n+1)$-dimensional compact manifold $X_0$. In this operation, we
see that the $n$-skeletal expansion $\bar{C}_\Delta$ of $\Delta$ is
exactly replaced by the connected sum between  $\bar{C}_\Delta$ and
$M'$. Thus, $N$ becomes actually a generalized $n$-skeletal
expansion of $(\Gamma, \alpha)$ such that $\Delta$ corresponds to
the interior of $M'\backslash B'$ rather than $C_\Delta$. It is not
difficult to see that the above argument still works for this $X_0$.
This completes the proof.
 \end{proof}

\begin{rem}
It should be pointed out that  if we start from {\em any} principal
$G$-bundle over $X_0$ (not just the trivial one), then  this
argument still works, i.e., we can always obtain a $G$-manifold
(c.f. \cite{J}, \cite{D}, \cite{D1}, and \cite{LM}).
\end{rem}

In the case where $n=2$, Theorem~\ref{realization} gives a complete
answer to question (Q4). This is because $(\Gamma, \alpha)$ in this
case always admits a 2-skeletal expansion, and the resulting surface
bounds a 3-manifold if and only if its Euler characteristic is even.

\begin{cor}
For every finite connected regular 3-valent graph $\Gamma$ with a
``good'' $G$-coloring $\alpha$ where $G=(\mathbb Z_2)^3$, there is a
3-dimensional $G$-manifold $X$ such that $(\Gamma,\alpha)$
corresponds to its moment graph.
\end{cor}

The case where $n=3$ is also very interesting because it is
well-known that every closed 3-manifold is the boundary of a
4-manifold.

\begin{cor}
For every finite connected regular 4-valent graph $\Gamma$ with a
``good'' $G$-coloring $\alpha$ where $G=(\mathbb Z_2)^4$, if
$\#\{3\text{-nests}\}=\#\{2\text{-nests}\}-\#\{\text{vertices}\}$,
then there is a 4-dimensional $G$-manifold $X$ such that
$(\Gamma,\alpha)$ corresponds to its moment graph.
\end{cor}


\begin{thebibliography}{GKM}
\bibitem[AP]{AP}{C. Allday, V. Puppe, \emph{Cohomological Methods in
 Transformation Groups}, Cambridge Studies in Advanced Mathematics, \textbf{32},
 Cambridge University Press, 1993.}
\bibitem[BGH]{BGH}{D. Biss, V. Guillemin, T. S. Holm, \emph{The mod $2$
 cohomology of fixed point sets of anti-symplectic involutions}, Adv. Math.
 \textbf{185} (2004), 370--399.}
\bibitem[CF]{CF}{P. E. Conner, E. E. Floyd, \emph{Differentiable periodic
 maps}, Springer-Verlag, 1964.}
\bibitem[D]{D}{M. Davis, \emph{Groups generated by reflections and aspherical
 manifolds not covered by Euclidean space}, Ann. of Math. \textbf{117} (1983),
 293--324.}
\bibitem[D1]{D1}{M. Davis, \emph{Smooth $G$-manifolds as collections of fiber
 bundles}, Pacific J. Math. \textbf{77} (1978), 315-363.}
\bibitem[DJ]{DJ}{M. Davis, T. Januszkiewicz, \emph{Convex polytopes, Coxeter
 orbifolds and torus actions}, Duke Math. J. \textbf{62} (1991), 417--451.}
\bibitem[GH]{GH}{V. Guillemin and T. Holm, \emph{GKM theory for torus actions
 with nonisolated fixed points}, Int. Math. Res. Not. \textbf{40} (2004),
 2105--2124.}
\bibitem[GHZ]{GHZ}{V. Guillemin, T. Holm and C. Zara, \emph{A GKM description
 of the equivariant cohomology ring of a homogeneous space}, J. Algebraic
 Combin. \textbf{23} (2006), 21--41.}
\bibitem[GKM]{GKM}{M. Goresky, R. Kottwitz, R. MacPherson, \emph{Equivariant
 cohomology, Koszul duality, and the localization theorem}, Invent. Math.
 \textbf{131} (1998), 25--83.}
\bibitem[GZ1]{GZ1}{V. Guillemin and C. Zara, \emph{Equivariant de Rham theory
 and graphs}, Asian J. Math. \textbf{3} (1999), 49--76.}
\bibitem[GZ2]{GZ2}{V. Guillemin, C. Zara, \emph{$1$-Skeleta, Betti numbers,
 and equivariant cohomology}, Duke Math. J. \textbf{107} (2001), 283--349.}
\bibitem[GZ3]{GZ3}{V. Guillemin and C. Zara, \emph{The existence of
 generating families for the cohomology ring of a graph}, Adv. Math.
 \textbf{174} (2003), 115--153.}
\bibitem[GZ4]{GZ4}{V. Guillemin, C. Zara, \emph{$G$-actions on graphs},
 Internat. Math. Res. Notices \textbf{10} (2001), 519--542.}
\bibitem[H]{H}{J. Hempel, \emph{$3$-Manifolds}, Annals Math. Studies
 \textbf{86}, Princeton U. Press, 1976.}
\bibitem[J]{J}{K. J\"anich, \emph{On the classification of $O(n)$-manifolds},
 Math. Ann. \textbf{176} (1968), 53--76.}
\bibitem[L]{L}{Z. L\"u, \emph{Graphs and $(\mathbb Z_2)^k$-actions},
 arXiv:math.AT/0508643.}
\bibitem[LM]{LM}{Z. L\"u, M. Masuda, \emph{Equivariant classification of
 $2$-torus manifolds}, preprint.}
\bibitem[MMP]{MMP}{H. Maeda, M. Masuda, T. Panov, \emph{Torus graphs and
 simplicial posets},  Adv. Math. \textbf{212} (2007), 458-483. arXiv: math.AT/0511582.}
\bibitem[MP]{MP}{M. Masuda, T. Panov, \emph{On the cohomology of torus
 manifolds}, Osaka J. Math. \textbf{43} (2006), 711--746.}
\bibitem[M]{M}{E. E. Moise, \emph{Geometric Topology in Dimensions $2$
 and $3$}, Springer, 1977.}
\bibitem[W]{W}{J. H. C. Whitehead, \emph{On $C^1$-complexes}, Ann. Math.
 \textbf{41} (1940), 809--824.}
\end{thebibliography}
\end{document}